\documentclass{birkjour_t2}

\pdfoutput=1

\usepackage[utf8]{inputenc}
\usepackage[english]{babel}

\usepackage{tabu}
\usepackage{makecell}

\usepackage{graphicx}

\usepackage{mathtools}
\usepackage{hyperref}
\usepackage{cleveref}

\usepackage{bm}
\newcommand{\vecc}[1]{\boldsymbol{\mathbf{#1}}}

\usepackage[backend=biber, sorting=none]{biblatex}
\title[ABM, DDEs and the Moon's Orbit]{%
On the Extension of Adams–Bashforth–Moulton Methods for Numerical
Integration of Delay Differential Equations
and Application to the Moon's Orbit
}

\author{Dan Aksim}
\address{Institute of Applied Astronomy of the Russian Academy of Sciences \\
    Russia, 191187, St. Petersburg, Kutuzova Embankment, 10}
\email{danaksim@iaaras.ru}

\author{Dmitry Pavlov}
\address{Institute of Applied Astronomy of the Russian Academy of Sciences \\
    Russia, 191187, St. Petersburg, Kutuzova Embankment, 10}
\email{dpavlov@iaaras.ru}

\begin{document}

\begin{abstract}
    One of the problems arising in modern celestial mechanics is the need of
    precise numerical integration of dynamical equations of motion of the
    Moon. The action of tidal forces is modeled with a time delay and the motion
    of the Moon is therefore described by a functional differential equation (FDE)
    called \emph{delay differential equation} (DDE).

    Numerical integration of the orbit is normally being performed in both
    directions (forwards and backwards in time) starting from some epoch (moment
    in time). While the theory of normal forwards-in-time numerical
    integration of DDEs is developed and well-known, integrating a DDE backwards
    in time is equivalent to solving a different kind of FDE called
    \emph{advanced differential equation}, where the derivative of the function
    depends on not yet known future states of the function.

    We examine a modification of Adams--Bashforth--Moulton method allowing to
    perform integration of the Moon's DDE forwards and backwards in time and the
    results of such integration.
\end{abstract}


\maketitle

\section{Introduction}

In 1960's, two things happened simultaneously: the space era began and computers
made their way into research labs. In years that followed, requirements to
accuracy of the coordinates of solar system bodies (ephemeris) for space
missions grew. Eventually, the so-called ``analytical'' ephemeris, calculated
via sophisticated series expansion, gave way to more precise numerical
ephemeris, obtained by numerical integration of dynamical equations on a
computer.

As the precision of astronomical observations improved, the dynamical models
underlying the ephemeris became more complex over time. The laser ranging of the
lunar retroreflectors, which began in the end of 1969, allowed to determine the
presence of a liquid core inside the Moon three decades later~\cite{williams01}.
Such findings require careful analysis of different causes of the observed
motion; specifically in the case of the Moon, it is important to model the
dissipation of energy that is caused by tides and the rotation of the elastic
Moon.

Redistribution of mass in a large rocky object does not happen instantly, but
rather takes hours to happen. Thus, the corresponding tidal terms in the
dynamical equations are naturally written with a time delay.

Let us move on to the equation for which the method descibed in this article was
developed. This is the equation for the rotational motion of the Moon
\cite{folkner,pavlov2016}, written in a Moon-fixed frame:
\begin{equation}
    \label{eq:euler}
    \dot{\vecc{\omega}} = \left( \frac{I}{m} \right)^{-1}
    \left[
        \frac{\vecc N}{m} - \frac{\dot I}{m} \vecc{\omega}
        - \vecc\omega\times\left(\frac{I}{m}\vecc\omega\right)
    \right],
\end{equation}
where $\vecc \omega$ is the angular velocity, $m$ is the mass of the Moon,
$\vecc N$ is the external torque, and $I(t)$ is the inertia tensor of the lunar mantle.
The inertia tensor is defined as
\begin{equation}
    \label{eq:inertia-tensor}
    \begin{aligned}
    \frac{I}{m} & = \frac{I_0}{m}
    - k_2\frac{\mu_\mathrm{E}}{\mu_\mathrm{M}}\left(\frac{R_\mathrm{M}}{r}\right)^5
    \left[\begin{array}{ccc}
    x^2-\frac{1}{3}r^2 & xy & xz \\
    xy & y^2-\frac{1}{3}r^2 & yz \\
    xz & yz & z^2-\frac{1}{3}r^2
    \end{array}\right] + \\
    & + k_2\frac{R_\mathrm{M}^5}{3\mu_\mathrm{M}}\left[\begin{array}{ccc}
    \omega_x^2-\frac{1}{3}(\omega^2-n^2) & \omega_x\omega_y &
    \omega_x\omega_z \\
    \omega_x\omega_y & \omega_y^2-\frac{1}{3}(\omega^2-n^2) &
    \omega_y\omega_z \\
    \omega_x\omega_z & \omega_y\omega_z & \omega_z^2 -
    \frac{1}{3}(\omega^2+2n^2) \\
    \end{array}\right],
    \end{aligned}
\end{equation}
where $I_0$ is the inertia tensor of the undistorted lunar mantle, 
$\mu_\mathrm{E}$ and $\mu_\mathrm{M}$ are the gravitational parameters
of Earth and the Moon respectively,
$R_\mathrm{M}$ is the equatorial radius of the Moon,
$k_2$ is the degree-2 Love number of the Moon,
$\bm r = (x, y, z)^T$ is the position of the Moon relative to Earth,
$n$ is the lunar mean motion.

Most importantly, the inertia tensor of the Moon is subject to
delayed tidal distortion from Earth and delayed spin distortion.
Tidal and spin distortions are evaluated with a delayed argument,
so calculation of $I(t)$ involves not $\bm r(t)$ and $\bm\omega(t)$,
but $\bm r(t-\tau)$ and $\bm\omega(t-\tau)$, and $\tau$ is a constant
time delay (in \cite{pavlov2016}, $\tau =  0.096\ \mathrm{d}$).
Thus, the equation \ref{eq:euler} is a
\emph{delay differential equation of neutral type with constant delays}
\cite{hale1993},
the general form for which would be
\begin{equation}
    \label{eq:dde_general}
    \dot{\vecc x}(t) = \vecc{f}(t,\, \vecc x(t),\, 
                                \vecc x(t - \tau),\, \dot{\vecc x}(t - \tau)),
\end{equation}
where $\vecc x =(x_{1}, x_{2},\dots , x_{n})$ and $\tau > 0$.

The initial condition for numerical integration of the Moon's orbit is of the form
\begin{equation}
    \label{eq:initial}
    \vecc{x}(t_0) = \vecc{x}_0
\end{equation}
and is defined for $t_0 = \mathrm{JD}\ 2446000.5$ (1984-10-27 midnight)
                   \cite{pavlov2016}.
To cover the whole interval of the lunar laser ranging observations,
numerical integration has to be performed \emph{both forward and backward in time}.

It is obvious that backward-in-time integration of DDE \eqref{eq:dde_general}
is equivalent to forward-in-time integration of
an \emph{advanced differential equation}
\begin{equation}
    \label{eq:ade_general}
    \dot{\vecc x}(t) = \vecc{f}(t,\, \vecc x(t),\, 
                                \vecc x(t + \tau),\, \dot{\vecc x}(t + \tau)).
\end{equation}

Therefore, what must be invented is a numerical method able to solve an initial value problem
\begin{equation}
    \begin{cases}
        \label{eq:dde_cauchy}
        \dot{\vecc x}(t) = \vecc{f}(t,\, \vecc x(t),\, 
                                \vecc x(t \pm \tau),\, \dot{\vecc x}(t \pm \tau)), \\
        \vecc{x}(t_0) = \vecc{x}_0.
    \end{cases}
\end{equation}
In other words, the method must be able to numerically integrate both
delay and advanced differential equations given an initial condition
\eqref{eq:initial}.

\section{Known approaches}

Most of existing numerical integrators for DDEs, like \cite{XPPAUT,DDE23}, work
only in forward direction in time, and require a so-called ``history function''
that unambiguously defines the $\vecc x(t-\tau)$ for $t\in[0,\tau)$. That
requirement conflicts with the nature of motion of celestial bodies,
which is usually quasiperiodic and does not ``take off'' from some initially
known trajectory.

Particularly for eq. \eqref{eq:euler}, some dedicated methods are known that
approximate the $\vecc x(t\pm\tau)$ and $\dot{\vecc x}(t\pm\tau)$.
In~\cite{Hofmann2018}, a quadratic approximation (based on known derivatives of
the equations) is used to calculate the delayed terms. In~\cite{folkner}, the
approximation is done via high degree polynomials that are pre-fit before the
main integration~\cite{WilliamsPersonal}.

\section{The ``nested integration'' method}

Suppose that for an arbitrary point $t_n = t_0 + n h$ the function value
$\vecc x_n = \vecc x(t_n)$ is known. Calculation of
$\vecc f_n = \vecc{f}(t_n,\, \vecc x_n,\,
             \vecc x(t_n \pm \tau),\, \dot{\vecc x}(t_n \pm \tau))$
then requires knowledge of delayed function and derivative values
$\vecc x(t_n \pm \tau)$ and $\dot{\vecc x}(t_n \pm \tau)$.

The ``nested integration'' method comes down to introducing
a new equation
\begin{equation}
    \label{eq:nodelay_equation}
    \dot{\vecc y}(t) = \vecc{g}(t,\, \vecc y(t)) = \vecc{f}(t,\, \vecc y(t),\, 
                                \vecc y(t),\, 0)
\end{equation}
and then finding an approximation to the delayed function value
$\vecc x(t_n \pm \tau)$ as
the numerical solution $\vecc y(t_n \pm \tau)$ of the initial value problem
involving eq.~\eqref{eq:nodelay_equation}
\begin{equation}
    \label{eq:nodelay_cauchy}
    \begin{cases}
        \dot{\vecc y}(t) = \vecc{g}(t,\, \vecc y(t)), \\
        \vecc y(t_n) = \vecc x_n.
    \end{cases}
\end{equation}
The delayed derivative $\dot{\vecc x}(t_n \pm \tau)$ in this case is simply
calculated as $\vecc{g}(t_n \pm \tau,\, \vecc y(t_n \pm \tau))$.
The value for $\dot{\vecc x}(t_n - \tau)$ in \cref{eq:nodelay_equation} is
chosen arbitrarily.

Both outer (of problem \eqref{eq:dde_cauchy}) and inner
(of problems \eqref{eq:nodelay_cauchy} occurring at each outer integration step)
integration can be performed by any numerical method, single- or multistep.

\section{The interpolation method}

Another way to calculate the delayed values of $\vecc x$ and $\dot{\vecc x}$
is to use polynomial interpolation of data computed on previous integration steps.
If the function values at nodes $[t_0, ..., t_p]$ are $[\vecc x_0, ..., \vecc x_p]$,
then for arbitrary point $t$ the value $\vecc x^{(p-1)}(t)$ of interpolating polynomial
of degree $(p-1)$ is given by the \emph{barycentric Lagrange formula}
\cite{berrut2004}
\begin{equation}
    \label{eq:lagrange}
    \vecc x^{(p-1)}(t) = \dfrac{\displaystyle \sum_{j=0}^p \dfrac{w_j}{t - t_j} \vecc x_j}
                   {\displaystyle \sum_{j=0}^p \dfrac{w_j}{t - t_j}},
\end{equation}
where $w_j$ are the \emph{barycentric weights} defined by
\begin{equation}
    \label{eq:lagrange_ws}
    w_j = \frac{1}{\prod_{k \neq j} (t_j - t_k)}, \ \ j = 0, ..., p.
\end{equation}
The weights $\{ w_j \}$ depend only on values of $\{ t_j \}$, or, for an equidistant grid
with $h = t_n - t_{n-1}$, on step size $h$.
Calculating the weights with \eqref{eq:lagrange_ws} requires $O(p^2)$
operations and has to be done once for one step size $h$.
Once the weights are known, evaluation of the Lagrange polynomial is performed
with \eqref{eq:lagrange} in $O(p \cdot n)$ operations.

\section{The ABMD Integrator}

The ABMD integrator, which is currently used in ERA-8 software system \cite{era8},
uses the Adams--Bashforth--Moulton predictor--corrector methods \cite{hairer_norsett}.

The Adams--Bashforth predictor and the Adams--Moulton corrector are given
by the formulas:
\begin{equation}
    \label{eq:adams_p}
    \vecc x_{n+1}^{\mathrm{(p)}} = \vecc x_n +
            h \sum_{j=0}^{k-1} \gamma_j \nabla^j \vecc f_n,
\end{equation}
\begin{equation}
    \label{eq:adams_c}
    \vecc x_{n+1} = \vecc x_n +
            h \sum_{j=0}^{k} c_j \nabla^j \vecc f_{n+1},
\end{equation}
where $k$ is the method order,
$\vecc x^{\mathrm{(p)}}$ is the predicted solution value,
and the coefficients $\gamma_j$ and $c_j$ are defined as
\begin{equation}
    c_j = - \sum_{i=0}^{j-1} \frac{c_i}{j + 1 - i},
\end{equation}
\begin{equation}
    \gamma_j = \sum_{k=0}^j c_k.
\end{equation}

The classic Adams--Moulton corrector formula \eqref{eq:adams_c} can be rewritten
in a more computationally efficient form.
First, we will extract $\vecc x_n$ from \eqref{eq:adams_p} and insert it
into \ref{eq:adams_c}:
\begin{equation}
    \vecc x_{n+1} = \vecc x_{n+1}^{\mathrm{(p)}} + h \left(
            \sum_{j=0}^{k} c_j \nabla^j \vecc f_{n+1} -
            \sum_{j=0}^{k-1} \gamma_j \nabla^j \vecc f_n
            \right).
\end{equation}
For convenience, let's denote
\begin{equation}
    \mathcal{S} = \sum_{j=0}^{k} c_j \nabla^j \vecc f_{n+1} -
                  \sum_{j=0}^{k-1} \gamma_j \nabla^j \vecc f_n.
\end{equation}
Now, substituting $c_j$ with $\gamma_j - \gamma_{j-1}$, taking out the last
term of the first sum and regrouping the sums:
\begin{equation}
    \begin{split}
        \mathcal{S} &=
        \sum_{j=0}^{k} (\gamma_j - \gamma_{j-1}) \nabla^j \vecc f_{n+1} -
                \sum_{j=0}^{k-1} \gamma_j \nabla^j \vecc f_n = \\
        &= \gamma_k \nabla^k \vecc f_{n+1} +
            \sum_{j=0}^{k-1} \gamma_j \left( \nabla^j \vecc f_{n+1} -
                                             \nabla^j \vecc f_n \right) -
            \sum_{j=0}^{k} \gamma_{j-1} \nabla^j \vecc f_{n+1}.
    \end{split}
\end{equation}
As $\nabla^j \vecc f_{n+1} - \nabla^j \vecc f_n = \nabla^{j+1} \vecc f_{n+1}$
and $\gamma_{-1} = 0$,
\begin{equation}
    \begin{split}
        \mathcal{S} &= \gamma_k \nabla^k \vecc f_{n+1} +
                    \sum_{j=0}^{k-1} \gamma_j \nabla^{j+1} \vecc f_{n+1} -
                    \sum_{j=1}^{k} \gamma_{j-1} \nabla^j \vecc f_{n+1} = \\
                &= \gamma_k \nabla^k \vecc f_{n+1} +
                \sum_{j=0}^{k-1} \gamma_j \nabla^{j+1} \vecc f_{n+1} -
                \sum_{i=0}^{k-1} \gamma_{i} \nabla^{i+1} \vecc f_{n+1} = \\
                &= \gamma_k \nabla^k \vecc f_{n+1}.
    \end{split}
\end{equation}
Hence, eq.~\ref{eq:adams_c} simplifies to
the \emph{modified Adams--Moulton corrector formula}
\begin{equation}
    \label{eq:adams_mc}
    \vecc x_{n+1} = \vecc x_n^{\mathrm{(p)}} + h \gamma_k \nabla^k \vecc f_{n+1},
\end{equation}
which, provided that the difference $\nabla^k \vecc f_{n+1}$ is known,
requires $O(n)$ operations instead of $O(k n)$ required by
the formula \eqref{eq:adams_c}.

So, the three-stage predictor--evaluator--corrector (PEC) algorithm comes down 
to calculating a predicted solution value $\vecc x_{n+1}^{\mathrm{(p)}}$
with eq.~\eqref{eq:adams_p}, evaluating the right-hand side
$\vecc f_{n+1} = \vecc f \left( t_{n+1},\ \vecc x_{n+1}^{\mathrm{(p)}} \right)$,
and then correcting the predicted solution value with eq.~\eqref{eq:adams_mc}.
The evaluator and corrector steps may be repeated multiple times:
such methods would be referred to as PECE, PECEC, PECECE etc.
The ABMD integrator uses a PECEC scheme.

By definition of the backward difference operator $\nabla$,
\begin{equation}
    \label{eq:backward_diff}
    \begin{aligned}
        \nabla \vecc f_i &= \vecc f_i - \vecc f_{i-1}, \\
        \nabla^k \vecc f_i &= \nabla^{k-1} \vecc f_i - \nabla^{k-1} \vecc f_{i-1}.
    \end{aligned}
\end{equation}
Therefore, calculating $k$ differences (from $\nabla \vecc f$
to $\nabla^k \vecc f$) requires $O(k^2 n)$ operations.
If, however, differences $\nabla \vecc f_i \dots \nabla^k \vecc f_i$
are known, then successively updating them to 
$\nabla \vecc f_{i+1} \dots \nabla^k \vecc f_{i+1}$
with the definition \eqref{eq:backward_diff} requires only $O(k n)$
operations, which makes it the preferred method of calculating backward
differences during integration, especially for problems with large values of $n$.

As follows from \cref{eq:adams_p,eq:adams_mc},
to compute a new solution value, a method of order $k$ requires
knowledge of $k$ previous RHS values.
At the start of the integration only one solution value is known, thus,
the other $k-1$ values have to be found with a startup procedure
using a single-step method.

In the ABMD integrator, startup procedure uses the nested integration method.
After the startup is finished, delayed solution values are retrieved with
Lagrange interpolation \eqref{eq:lagrange}, which turns to \emph{extrapolation}
when integration is performed backward in time.

\section{Results}

We did not find DDEs that are sensible and have known solutions or invariants to
check against. So the primary criterion to validate the proposed method will be
the numerical integration of the ``real world''
\cref{eq:euler,eq:inertia-tensor}.

We used the ABMD integrator with Adams--Bashforth--Moulton PECEC scheme of order
13 with step size and with step size $h_\mathrm{ABM} = 1/16 \ \mathrm{d}$.

The startup uses the 8th order Dormand--Prince method~\cite{hairer_norsett}
for the outer integration method and the 4th order Runge--Kutta for the inner
integration. To match the precision of the 13th order ABM, step size
for Dormand--Prince startup is reduced: $h_\mathrm{DP} = h_\mathrm{ABM} / 8$.

\begin{figure}
    \includegraphics{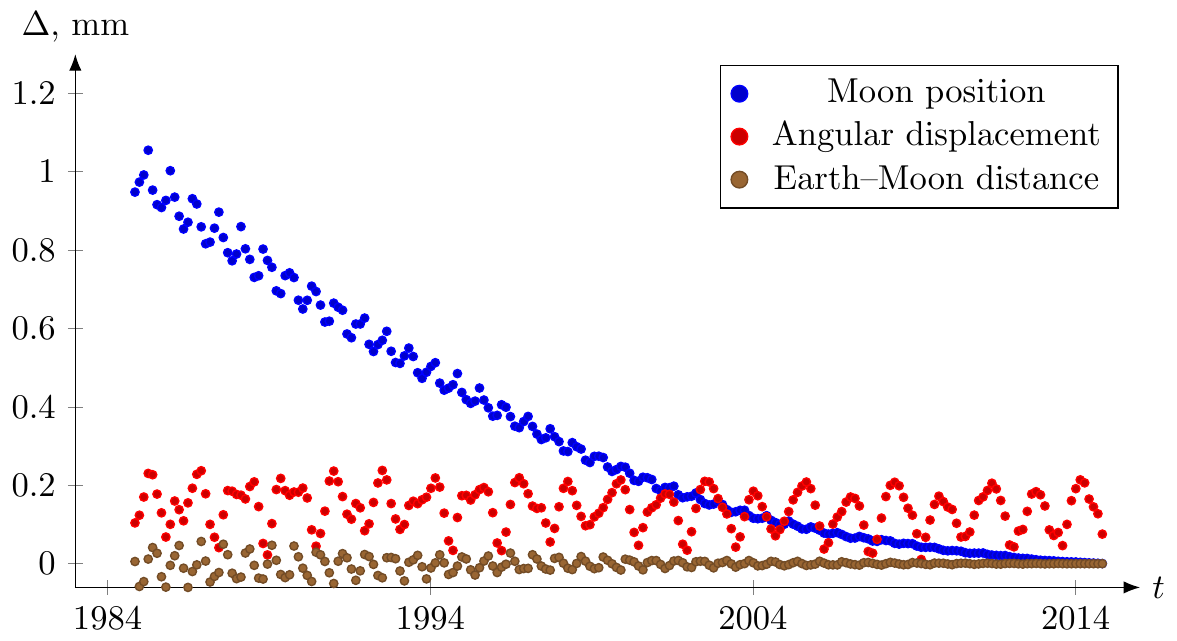}
    \caption{Differences between Moon positions, Earth--Moon distances and
      Moon's rotations (in maximum displacement on surface) acquired with
      forward and backward integrations.}
    \label{fig:forward-backward}
\end{figure}

Two tests were done. The first one is the usual forward--backward test.  The
full dynamical equations of the Moon from~\cite{pavlov2016} were integrated from
1984 to 2014, and the solution at the last point was taken as the initial
condition for integrating backwards, from 2014 to 1984. Differences between the
two results are shown in \cref{fig:forward-backward}. It is evident that the
inaccuracy inherent to the method is sufficient at present level of best lunar
laser ranging observations (2 mm in Earth--Moon distance).

\begin{figure}
    \centering
    \includegraphics{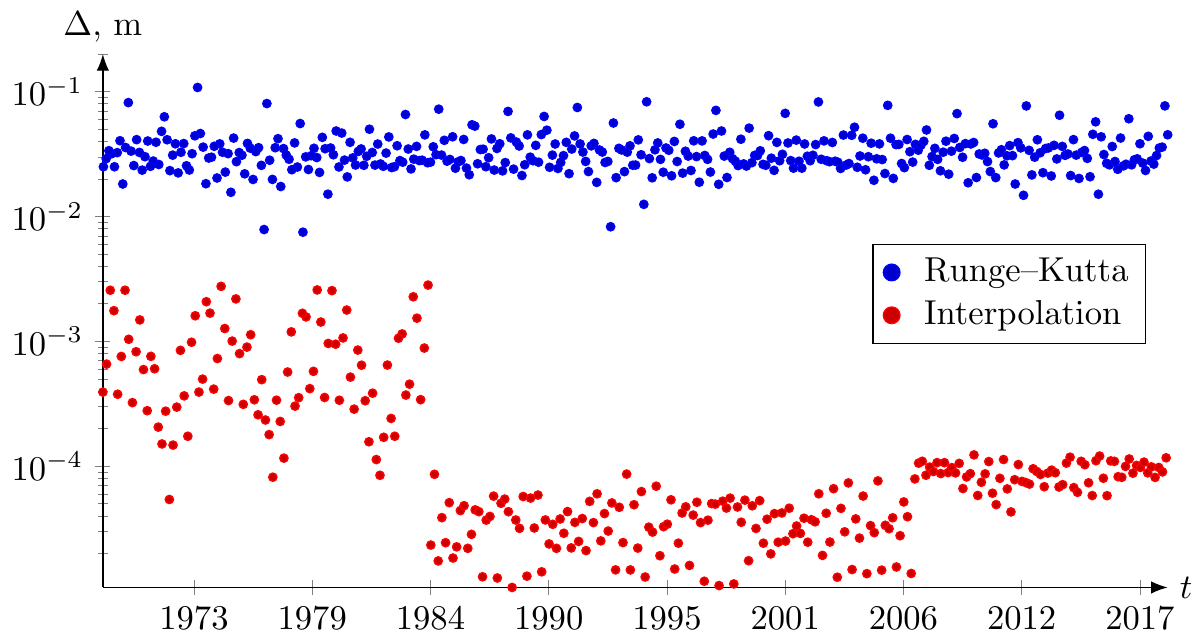}
    \caption{Differences between the delayed values obtained during integration
             and those of the final orbit.}
    \label{fig:delayed-diff}
\end{figure}

Another test was done to compare the delayed values $\vecc r(t - \tau)$ that are obtained
during integration by interpolation with
the same values of the final integrated orbit. Fig.~\ref{fig:delayed-diff} shows
that values calculated with polynomial interpolation---forward from the epoch---are
two orders of magnitude closer to the final orbit than those that are calculated
with nested Runge--Kutta integration. For the values calculated with
extrapolation---backward from the epoch---the improvement is one order.

\section{Conclusion and further work}

A modification of the Adams--Bashforth--Moulton multistep numerical integration
method was proposed, and its implementation, ABMD, was developed that allows to
integrate delay differential equations.  Unlike most other methods for DDEs, the
developed method does not require a ``history function'' and can integrate
forwards and backwards in time.

There has been no theoretical analysis of the method's correction or
convergence. Most probably the applicability of ABMD is limited to relatively
stable systems in which the delay term is not the main term of the equation.  It
has been shown that ABMD works well for one particular equation of the
Earth--Moon system.

The proposed method is trivially generalised to DDEs with multiple constant
delays with possibly different signs. The present model of the Earth--Moon
system already includes delays in other equations besides the one that was
described in this paper: these are equations for the Earth
tides~\cite{Williams2016}.

For spacecraft dynamics (as opposed to dynamics of celestial bodies), more
precise and complicated tidal models are used, see e.g.~\cite{WilliamsBoggs2015}
for lunar tides. Such models do not contain the delays directly in the
equations, but instead use empirical Fourier series whose arguments are linear
combinations of Delaunay variables linked to the state of the Earth--Moon--Sun
system.

We believe that the new tool developed for integration of arbitrary DDEs can
help to build more natural models for tides on Earth and the Moon. Such new
models could be a step to a more adequate description of the mechanical
processes that are going at the interior of the celestial bodies.

\section{Acknowledgements}

Authors are grateful to Sergey Kurdubov (IAA RAS) and John Chandler (Harvard
CfA) for useful discussions related to this work. Steve Moshier's implementation
of Adams--Bashforth--Moulton method (\url{http://moshier.net/ssystem.html}) was
a helpful reference.

\printbibliography

\end{document}